\documentclass [11pt]{article}

\usepackage[reqno]{amsmath}
\usepackage{amsfonts}
\usepackage{amssymb}
\usepackage{amsthm}
\usepackage{amscd}
\usepackage{graphicx,color}
\usepackage{newlfont,color}
\usepackage{dsfont}
\usepackage{pifont}
\usepackage{epstopdf}
\usepackage{natbib}

\newcommand{\p}{{\mathbb P}}
\newcommand{\e}{{\mathbb E}}

\newcommand{\eqb}{\begin{equation}}
\newcommand{\eqe}{\end{equation}}

\newcommand{\Proof}{{\bf \noindent Proof: }}

\theoremstyle{definition}          \newtheorem{defn}{Definition}[section]
\theoremstyle{definition}          \newtheorem{prob}{Problem}[section]
\theoremstyle{remark}              \newtheorem{rem}{Remark}[section]
\theoremstyle{plain}               \newtheorem{thm}{Theorem}[section]
\theoremstyle{plain}               
\theoremstyle{plain}               \newtheorem{pro}{Proposition}[section]
\theoremstyle{plain}               
\theoremstyle{definition}          
\theoremstyle{remark}              

\begin{document}
\pagenumbering{arabic}             

\pagestyle{plain}


\title{A maximum principle for Markov-modulated SDEs of mean-field type and mean-field game}

\author{Yongming Tai}

\date{Department of Management Sciences \\University of Waterloo
       \\ Waterloo, Ontario, Canada N2L 3G1
       \\ May 19, 2014}

\maketitle                        

\begin{abstract}
In this paper, we analyze mean-field game modulated by finite states markov chains.
We first develop a sufficient stochastic maximum principle for the optimal control of a Markov-modulated
stochastic differential equation (SDE) of mean-field type whose coefficients depend on
the state of the process, some functional of its law as well as variation of time and sample.
As coefficients are perturbed by a Markov chain and thus random, to study such SDEs, we analyze
existence and uniqueness of solutions of a class of mean-field type SDEs whose coefficients are random
Lipschitz as well as the property of propagation of chaos for associated interacting particles system
with method parallel to existing results as a byproduct. We also solve approximate Nash equilibrium
for the Markov-modulated mean-field game by mean-field theory.

\vspace{2mm}

\noindent $Keywords:$ Mean-field type SDEs, interacting particle systems, Markov-modulated, stochastic control,
stochastic maximum principle, mean-field game, approximate Nash equilibrium

\end{abstract}

\section{Introduction}\label{mpmm2014int}

Mean-field type stochastic differential control (games) has become a very popular topic and developed rapidly
during recent years, for instance, see \citet{huang2003individual,huang2006large} and \citet{lasry2007mean}.
Tools to analyze such problems include mean-field games method \citep[e.g. see][]{huang2006large}
which solves a standard control
problem for  a deterministic function firstly, and then utilize fixed point method to determine
that there exists such function which is the distribution of the state process, and partial differential
equation~\citep[e.g. see][]{borkar2010mckean} and stochastic maximum principle of mean-field type
\citep[e.g. see][]{andersson2011maximum} which solve mean-field type control problem in a direct manner.

This paper aims to generalize
the stochastic maximum principle of mean-field type developed in \citet{andersson2011maximum}
to the case modulated by a homogenous Markov chain defined finite state space, and analyze
Markov-modulated mean-field game. We consider the stochastic problem of a
Markov-modulated SDE of mean-field type whose coefficients depend on
the state of the process, some functional of its law as well as variation of time and sample (see
Problem~\ref{mpmm2014prob1} for details). Under suitable assumptions, the Markov-modulated SDE of mean-field type
can be obtained as a limit of an interacting particles system modulated by independent, identically
distributed (i.i.d.) Markov chains defined on finite state space (see the system~(\ref{mpmm2014eqn1})
for details). We provide sufficient conditions for maximum principle of mean-field type modulated by
a Markov chain. Compared to \citet{andersson2011maximum}, the adjoint equation in this paper
involves Markov regime-switching jumps which
inhabit property of martingale. As in \citet{andersson2011maximum}, we establish the sufficient conditions
for maximum principle in general set-up that only requires predictable processes and some
integrability conditions without any other special conditions for admissible controls, which makes
it suitable to handle time inconsistency of the mean-field type control problem. Stochastic maximum
principles of mean-field type include \citet{buckdahn2011general}, \citet{hosking2012stochastic},
\citet{li2012stochastic}, \citet{shen2013maximum}. For standard
maximum principles modulated by Markov chains, please see \citet[e.g.][]{zhang2012stochastic}.

As coefficients of mean-field type SDEs are left continuous with right limit, and
perturbed by a Markov chain and thus random,
to study such SDEs, it is necessary to analyze
existence and uniqueness of solutions of a class of mean-field type SDEs whose coefficients are random
Lipschitz as well as the property of propagation of chaos for associated interacting particles system
although methods used by us are simply parallel to the ones from existing results in \citet{jourdain2007nonlinear}.

We also analyze asymptotic errors between stochastic dynamic game and its limiting optimal problem.
By mean-field game theory, we find approximate Nash equilibrium for our model.

The remaining of this paper is organized as follows. In Section~\ref{mpmm2014prof}, we introduce
Markov-modulated diffusion, weakly coupled stochastic dynamic game and mean-field type SDEs,
and analyze existence and uniqueness of mean-field type SDEs as well as the property of
propagation of chaos of the associated interacting particles system. In Section~\ref{mpmm2014scmp},
we provide sufficient conditions for maximum principle. In Section~\ref{mpmm2014ane}, we make
make asymptotic analysis and verify approximate Nash equilibrium. In Section~\ref{mpmm2014conr}, we
comment concluding remarks.

\section{Problem formulation}\label{mpmm2014prof}

In this section, we define the Markov-modulated diffusion model and formulate the
Markov-modulated weakly coupled interacting particles control system governed by $n$ dimensional
nonlinear stochastic system and its limit system governed by Markov-modulated mean-field type
SDEs. Meanwhile, we analyze existence and uniqueness of two systems as well as approximation of
two systems. Assume that $T>0$ is time horizon and $(\Omega, \mathcal{F}, \p)$ is a complete
probability space. We shall make use of the following notation:\\
\indent $\mathds{R}$\hspace{20mm}:\quad the real Euclidean space;\\
\indent $\mathds{R}^{n}$\hspace{18mm}:\quad the $n$-dimensional real Euclidean space;\\
\indent $|\cdot|$\hspace{18mm}:\quad the Euclidean norm;\\
\indent $M^{*}$\hspace{17mm}:\quad the transpose of any matrix or vector;\\
\indent $Diag(y)$\hspace{9mm}:\quad the diagonal matrix with the elements of $y$;\\
\indent $L_{\mathcal{M}}^{2}(0, T; \mathds{V})$\hspace{1mm} :\quad the space of all measurable,
                        $\mathcal{M}_{t}$-predictable processes\\
          \indent\hspace{27mm} $f: [0,T]\times\Omega\longrightarrow \mathds{V}$ such that
                                  $\e\int_{0}^{T}|f(t)|^{2}dt < \infty$, where\\
          \indent\hspace{27mm} $\mathcal{M}_{t}$ is a $\sigma$-field and $\mathds{V}$ is a subset of $\mathds{R}$;\\
\indent $\mathds{U}$\hspace{20mm}:\quad the action space which is a nonempty, closed and convex
               subset of $\mathds{R}$;\\
\indent $C$\hspace{20mm}:\quad a constant which may change from line to line.

\subsection{The Markov-modulated diffusion model}

Let the time-homogeneous Markov chain $\alpha$ take values in finite state space
$\mathds{S}\triangleq\{1,2,\dots,d\}$ associated with the generator
$\Lambda\triangleq \left[\lambda_{ij}\right]_{i,j=1,\dots,d}$. We assume that that chain starts in a fixed state
$i_{0}\in\mathds{S}$ such that $\alpha_{0}=i_{0}$. Let $N_{t}(i,j)$ be the number of jumps from
state $i$ to state $j$ up to time $t$. Denote by $\boldsymbol{1}$ the indicator function. Then
\[
N_{t}(i,j) = \sum\limits_{0<s\leq t}\boldsymbol{1}_{\{\alpha_{s-}=i\}}\boldsymbol{1}_{\{\alpha_{s}=j\}},\,\,
\forall t\in [0,T].
\]
Define the intensity process $m_{t}(i,j)\triangleq \lambda_{ij}\boldsymbol{1}_{\{\alpha_{t-}=i\}}$.
If we compensate $N_{t}(i,j)$ by $\int_{0}^{t}m_{s}(i,j)ds$, then the resulting process
$N_{t}(i,j)-\int_{0}^{t}m_{s}(i,j)ds$ is a purely discontinuous, square integrable martingale
which is null at the origin. Let
\[
\widetilde{\Phi}_{t}(j)\triangleq \sum\limits_{i=1, i\neq j}^{d}\left(N_{t}(i,j)-\int_{0}^{t}m_{s}(i,j)ds\right),\,\,
j=1,\dots,d,
\]
$\widetilde{\Phi}_{t}(j)$ is a martingale. Denote $\widetilde{\Phi}_{t}\triangleq
(\widetilde{\Phi}_{t}(1),\dots,\widetilde{\Phi}_{t}(d))$ and $m_{t}\triangleq(m_{t}(1),\dots,m_{t}(d))$,
where $m_{t}(j)\triangleq \sum\limits_{i=1, i\neq j}^{d}\int_{0}^{t}m_{s}(i,j)ds$, for $j=1,\dots,d$.

\subsection{Weakly coupled stochastic dynamic game and mean-field type SDEs}

We consider weakly coupled system of $n$ interacting particles modulated by independent Markov
chains. The dynamic of each particle is given by
\eqb\label{mpmm2014eqn1}
\left\{
  \begin{array}{lll}
    dx_{t}^{i,n} & = & b\left(t, x_{t-}^{i,n}, \frac{1}{n}\sum\limits_{j=1}^{n}\psi(x_{t-}^{j,n}), u_{t}^{i}\right)
                              r\left(\alpha_{t-}^{i}\right)dt
                          +\sigma\left(t, x_{t-}^{i,n}, \frac{1}{n}\sum\limits_{j=1}^{n}\phi(x_{t-}^{j,n}), u_{t}^{i}\right)
                              r\left(\alpha_{t-}^{i}\right)dw_{t}^{i}, \\
    x_{0}^{i,n}  & = & x_{0}^{i},\,\,\, i=1,\dots,n,\quad t\in [0,T],
  \end{array}
\right.
\eqe
where $w_{t}^{1},\dots,w_{t}^{n}$, $t\in [0,T]$, are $n$ independent standard scalar Brownian
motions and $\alpha_{t}^{1},\dots,\alpha_{t}^{n}$, $t\in [0,T]$, are $n$
i.i.d. time-homogeneous Markov chains taking values in $\mathds{S}$.
The initial states $x_{0}^{1},\dots,x_{0}^{n}$ are mutually
independent and satisfy $\e|x_{0}^{i}|^{2}<\infty$, $i=1,\dots,n$. We also assume that
$\{x_{0}^{1},\dots,x_{0}^{n}\}$, $\{w_{t}^{1},\dots,w_{t}^{n}\}$ and
$\{\alpha_{t}^{1},\dots,\alpha_{t}^{n}\}$ are mutually independent. The control
$u^{i}\in L_{\mathcal{G}}^{2}(0, T; \mathds{U})$, $i=1,\dots,n$, where
\[
\mathcal{G}_{t} \triangleq
\sigma(\{x_{0}^{1},\dots,x_{0}^{n}\}, \{w_{s}^{1}, \dots, w_{s}^{n}\},
\{\alpha_{t}^{1},\dots,\alpha_{t}^{n}\}, s \leq t).
\]
The functions $b$, $\sigma$, $\psi$ and $\phi$ are given as follows:
\[
\begin{array}{lll}
  b      & : & [0,T]\times\mathds{R}\times\mathds{R}\times\mathds{U}\longrightarrow\mathds{R}, \\
  \sigma & : & [0,T]\times\mathds{R}\times\mathds{R}\times\mathds{U}\longrightarrow\mathds{R},  \\
  \psi   & : & \mathds{R}\longrightarrow\mathds{R}, \\
  \phi   & : & \mathds{R}\longrightarrow\mathds{R},\\
  r      & : & \mathds{R}\longrightarrow\mathds{R}.
\end{array}
\]
Each particle is often called an agent (or a player).

The cost functional for the $i$th agent is given by
\[
\begin{array}{ll}
             &\mathcal{J}^{i}(u^{1}, \dots, u^{n}) \\
\triangleq  &
\e\left(\int_{0}^{T}h\left(t, x_{t}^{i,n}, \frac{1}{n}\sum\limits_{j=1}^{n}\varphi(x_{t}^{j,n}),
u_{t}^{i}\right)r\left(\alpha_{t}^{i}\right)dt
 + g\left(x_{T}^{i,n}, \frac{1}{n}\sum\limits_{j=1}^{n}\chi(x_{T}^{j,n})\right)r\left(\alpha_{T}^{i}\right)\right),
\end{array}
\]
where
\[
\begin{array}{lll}
  g      & : & \mathds{R}\times\mathds{R}\longrightarrow\mathds{R}, \\
  h      & : & [0,T]\times\mathds{R}\times\mathds{R}\times\mathds{U}\longrightarrow\mathds{R},  \\
  \varphi   & : & \mathds{R}\longrightarrow\mathds{R}, \\
  \chi   & : & \mathds{R}\longrightarrow\mathds{R}.
\end{array}
\]

The objective of each agent is to minimize his own cost by properly controlling his own dynamics.
Due to the interaction between agents, the computation of a Nash equilibrium is highly complicated,
especially for a large population of agents. In practice, a convenient computable strategy is highly
demanded. Instead of Nash equilibrium, an approximate $\varepsilon$-Nash equilibrium which was
introduced successfully solve this problem \citep[e.g. see][]{huang2006large}.
\begin{defn}\label{mpmm2014defn1}
For the $n$ agents, a sequence of controls $u^{i}\in L_{\mathcal{G}}^{2}(0, T; \mathds{U})$
$(resp.,u^{i}\in L_{\mathcal{F}^{i}}^{2}(0, T; \mathds{U}))$ which is a Lipschitz feedback, where
$\mathcal{F}^{i}_{t} \triangleq \sigma(x_{0}^{i}, w_{s}^{i}, \alpha_{t}^{i}, s \leq t)$,
$i=1,\dots,d$, is called
$\varepsilon$-Nash equilibrium with respect to the cost $\mathcal{J}^{i}(u^{1}, \dots, u^{n})$
if there exists $\varepsilon>0$ such that for any fixed $1\leq i \leq n$, we have
\[
\mathcal{J}^{i}(u^{1}, \dots, u^{n}) \leq \mathcal{J}^{i}(u^{1}, \dots, u^{i-1}, v^{i}, u^{i+1}, \dots, u^{n})
+\varepsilon,
\]
when any alternative control $v^{i}\in L_{\mathcal{G}}^{2}(0, T; \mathds{U})$
$(resp.,v^{i}\in L_{\mathcal{F}^{i}}^{2}(0, T; \mathds{U}))$ which is another Lipschitz feedback
is applied by the $i$th agent.
\end{defn}
%
By the mean field game theory, a candidate for $\varepsilon$-Nash equilibrium can
be solved via solving the following limiting problem.
\begin{prob}\label{mpmm2014prob1}
Find an control strategy $\hat{u}\in L_{\mathcal{\bar{F}}}^{2}(0,T;\mathds{U})$,
minimize
\[
             J(\bar{u}) \triangleq
\e\left(\int_{0}^{T}h\left(t, x_{t}, \e\varphi(x_{t}),
\bar{u}_{t}\right)r\left(\alpha_{t}\right)dt
 + g\left(x_{T}, \e\chi(x_{T})\right)r\left(\alpha_{T}\right)\right)
\]
subject to
\eqb\label{mpmm2014eqn2}
\left\{
  \begin{array}{lll}
    dx_{t} & = & b(t, x_{t-}, \e\psi(x_{t-}), \bar{u}_{t})
                              r(\alpha_{t-})dt
                          +\sigma(t, x_{t-}, \e\phi(x_{t-}), \bar{u}_{t})
                              r(\alpha_{t-})dw_{t}, \\
    x_{0}  & = & x(0),
  \end{array}
\right.
\eqe
for any control $\bar{u}\in L_{\mathcal{\bar{F}}}^{2}(0,T;\mathds{U})$, where
$\bar{\mathcal{F}}_{t} \triangleq \sigma(x(0), w_{s}, \alpha_{t}, s \leq t)$.
We assume that $w_{t}, t\in [0,T]$, is a standard scalar Brownian motion. $\alpha_{t}, t\in [0,T]$, is a
time-homogeneous Markov chain defined on $\mathds{S}$ and independent of $w_{t}$. The initial state
satisfies $\e|x(0)|^{2}<\infty$ and is independent of $w_{t}$ and $\alpha_{t}$.
In the above equation, the expectation means the conditional expectation conditioned on
$\{x_{0}=x(0), \alpha_{0}=i_{0}\}$.
\end{prob}

Once Problem~\ref{mpmm2014prob1} were solved
with $\hat{u}$, we could obtain controls $u^{i}\in L_{\mathcal{F}^{i}}^{2}(0,T;\mathds{U})$,
$i=1,\dots,n$.
$u^{1},\dots,u^{n}$ are independent. It can be shown that $(u^{1},\dots,u^{n})$
is a $\varepsilon$-Nash equilibrium. In this paper, we shall develop a stochastic maximum
principle of mean-field type for Problem~\ref{mpmm2014prob1} which generalizes the result
in \cite{andersson2011maximum}. The following
assumptions will be imposed throughout this paper, where $x$ denotes the state variable, $y$
the "expected value", $v$ the control variable and t the time.
\begin{description}
\item[(A.1)] $\psi$, $\phi$, $\chi$ and $\varphi$ are continuously differentiable. $g$ is
             continuously differentiable with respect to $(x,y)$. $b$, $\sigma$ and $h$ are
             continuously differentiable with respect to $(x,y,v)$. $b$ and $\sigma$ are
             left continuous with right limit with respect to $t$. $r$ is positive continuous
             function.
\item[(A.2)] All the derivatives in (A.1) are Lipschitz continuous and bounded.
\item[(A.3)] $\int_{0}^{T}(|b(s,0,\psi(0),0)|^{2} + |\sigma(s,0,\phi(0),0)|^{2})ds < \infty$.
\end{description}
Under the assumptions (A.1), (A.2) and (A.3), the solutions to the system~(\ref{mpmm2014eqn1}) and equation~(\ref{mpmm2014eqn2}) are unique.
\begin{rem}
If $b$ and $\sigma$ are continuous with respect to $t$, then the assumption (A.3) can be relaxed.
\end{rem}

\subsection{Solution to mean-field type SDEs and propagation of chaos}

In this subsection, we analyze solvability of the system~(\ref{mpmm2014eqn1}) and the equation
(\ref{mpmm2014eqn2}). Instead of directly analyzing the
system~(\ref{mpmm2014eqn1}) and the equation~(\ref{mpmm2014eqn2}), we study a more general case
driven by square integrable L\'{e}vy processes $\{z_{t}, t\in [0,T]\}$ parallel to
\citet{jourdain2007nonlinear}. To do this, we introduce nonlinear
stochastic differential equation of mean-field type modulated by a Markov chain defined on
$\mathds{S}$ and corresponding system of $n$ interacting particles as follows:
\eqb\label{mpmm2014eqn3}
\left\{
  \begin{array}{lll}
    d\tilde{x}_{t} & = & \tilde{\sigma}(t, \tilde{x}_{t-}, \p_{t-}))
                              r(\tilde{\alpha}_{t-})dz_{t}, \quad t\in [0,T]\\
    \tilde{x}_{0}  & = & \tilde{x}(0),
  \end{array}
\right.
\eqe
where $\{z_{t}, t\in [0,T]\}$ is a L\'{e}vy process with value in $\mathds{R}$,
$\{\tilde{\alpha}_{t}, t\in [0,T]\}$ is a Markov chain on $\mathds{S}$.
For $t\in [0,T],\, \p_{t}\,\, \text{denotes the probability
                        distribution of}\,\, \tilde{x}_{t}$, and
$\mathbb{P}_{s-}=\mathbb{P}\circ\tilde{x}_{s-}^{-1}$ is the weak limit of $\mathbb{P}_{t}$
as $t\to s$ increasingly. The initial state
$\tilde{x}(0)$ takes values in $\mathds{R}$, distributed according to $\pi$,
and satisfies $\e|\tilde{x}(0)|^{2} < \infty$.
Furthermore, $\tilde{x}(0)$, $\{z_{t}, t\in [0,T]\}$ and $\{\tilde{\alpha}_{t}, t\in [0,T]\}$ are mutually
independent. The functions $\tilde{\sigma}$ is given as follows:
\[
  \tilde{\sigma} : [0,T]\times\mathds{R}\times\mathcal{P}(\mathds{R})\longrightarrow\mathds{R},
\]
where $\mathcal{P}(\mathds{R})$ is the space of probability measures on $\mathds{R}$. By choosing
$\tilde{\sigma}$ linear in the third variable, the classical Mckean-Vlasov model studied in
\citet{sznitman1991topics} can be obtained as a special case of~(\ref{mpmm2014eqn3}). Let
$\tilde{\tilde{\sigma}}(t,\omega,\cdot,\p_{t-})\triangleq\tilde{\sigma}(t, \cdot, \p_{t-})r(\tilde{\alpha}_{t-}(\omega))$.
It is noted that the above equation
have more general coefficient which depends not only on the state process and the probability distribution
of the state process, but also on the time and the sample.

For $i=1,\dots,n$, let
$(\tilde{x}_{0}^{i}, z^{i})$ be a sequence of independent copies of $(\tilde{x}(0), z)$.
Define the weakly coupled system
of $n$ interacting particles
\eqb\label{mpmm2014eqn4}
\left\{
  \begin{array}{lll}
    d\tilde{x}_{t}^{i,n} & = & \tilde{\sigma}(t, \tilde{x}_{t-}^{i,n}, \mu_{s-}^{n}))
                              r(\tilde{\alpha}_{t-})dz_{t}^{i}, \quad t\in [0,T],\,\, i=1,\dots,n,\\
    \tilde{x}_{0}^{i,n}  & = & \tilde{x}_{0}^{i},\,\,\, \mu^{n}\triangleq
    \frac{1}{n}\sum_{j=1}^{n}\delta_{\tilde{x}^{j,n}}\,\, \text{is the empirical distribution},\,\,
    \delta_{x}\,\,\text{is Dirac measure}.
  \end{array}
\right.
\eqe
We shall show that as $n\to\infty$, for $i=1,\dots,n$, $\tilde{x}^{i,n}$ converge to a limit
$\tilde{x}^{i}$ which is an independent copy of the solution to equation~(\ref{mpmm2014eqn3}).

Let $\mathcal{P}_{2}(\mathds{R})$ be the space of probability measures on $\mathds{R}$ with finite second order
moments. For $\mu, \nu\in\mathcal{P}_{2}(\mathds{R})$, we define the Vaserstein metric as follows:
\[
d(\mu,\nu)=\inf\left\{\left(\int_{\mathds{R}\times\mathds{R}}|x-y|^{2}Q(dx,dy)\right)^{1/2}: Q\in
\mathcal{P}(\mathds{R}\times\mathds{R})\,\, \text{with marginals}\,\, \mu\,\, \text{and}\,\, \nu\right\}.
\]
It induces the topology of weak convergence together with convergence of moments up to order 2.
Due to $r(\cdot)$ being continuous and thus bounded on any compact set, if for each $t$,
$\tilde{\sigma}(t,\cdot,\cdot)$ is Lipschitz continuous when $\mathds{R}\times\mathcal{P}_{2}(\mathds{R})$
is endowed with the product of the canonical metric on $\mathds{R}$ and the Vaserstein metric
on $\mathcal{P}_{2}(\mathds{R})$, then for each $(t,\omega,\cdot,\nu)$, $\tilde{\tilde{\sigma}}(t,\omega,\cdot,\nu)$
is random Lipschitz continuous with respect to the canonical metric on $\mathds{R}$.

The solvability of the equation~(\ref{mpmm2014eqn3}) and the system~(\ref{mpmm2014eqn4}) is
given as follows:
\begin{pro}\label{mpmm2014pro1}
Assume that $\{z_{t}, t\in [0,T]\}$ is square integrable, and that for each $t$,
$\tilde{\sigma}(t,\cdot,\cdot)$ is Lipschitz continuous when $\mathds{R}\times\mathcal{P}_{2}(\mathds{R})$
is endowed with the product of the canonical metric on $\mathds{R}$ and the Vaserstein metric
on $\mathcal{P}_{2}(\mathds{R})$. In addition, for fixed $x$ and $\nu$, $\tilde{\sigma}(\cdot, x,\nu)$
is left continuous with right limit and $\int_{0}^{T}|\tilde{\sigma}(s,0,\delta_{0})|^{2}ds<\infty$, where
$\nu\in\mathcal{P}_{2}(\mathds{R})$ and $\delta_{0}$ is Dirac measure.
Then equation~(\ref{mpmm2014eqn3}) admits a unique strong solution such that
$\e(\sup_{t\leq T}|\tilde{x}_{t}|^{2})<\infty$.
\end{pro}
\noindent Since for
$\xi=(x_{1},\dots,x_{n})$ and $\zeta=(y_{1},\dots,y_{n})$ in $\mathds{R}^{n}$, we have
\eqb\label{mpmm2014eqn5}
d\left(\frac{1}{n}\sum_{j=1}^{n}\delta_{x_{j}},\frac{1}{n}\sum_{j=1}^{n}\delta_{y_{j}}\right)
\leq \left(\frac{1}{n}\sum_{j=1}^{n}|x_{j}-y_{j}|^{2}\right)^{1/2}=\frac{1}{\sqrt{n}}|\xi-\zeta|.
\eqe
Thus, let $\vartheta(t,\omega,x_{1},\dots,x_{n})\triangleq
\tilde{\tilde{\sigma}}(t,\omega,x_{i},\frac{1}{n}\sum_{j=1}^{n}\delta_{x_{j}})$, we have
$\vartheta: [0,T]\times\Omega\times\mathds{R}^{n}\longrightarrow\mathds{R}$ is random
Lipschitz continuous. It induces a process Lipschitz operator which is therefore functional
Lipschitz. Hence, existence of a unique solution to the system~(\ref{mpmm2014eqn4}), with
finite second order moments, follows from Theorem 7, p.253, in \citet{protter2004stochastic}.

Next, we give the trajectorial propagation of chaos result for the system~(\ref{mpmm2014eqn4}).
\begin{pro}\label{mpmm2014pro2}
Under the assumptions of Proposition~\ref{mpmm2014pro1}
\[
\lim_{n\to\infty}\sup_{i\leq n}\e\left(\sup_{t\leq T}|\tilde{x}_{t}^{i,n}-\tilde{x}_{t}^{i}|^{2}\right)=0.
\]
Moreover, if $\tilde{\sigma}(t,\tilde{x},\nu)=\int_{\mathds{R}}\eta(t,\tilde{x},\tilde{y})\nu(d\tilde{y})$, where
$\eta: [0,T]\times\mathds{R}\times\mathds{R}\longrightarrow\mathds{R}$ is a Lipschitz continuous function
with respect to $(\tilde{x},\tilde{y})$ and left continuous with right limit with respect to $t$,
then
\[
\sup_{i\leq n}\e\left(\sup_{t\leq T}|\tilde{x}_{t}^{i,n}-\tilde{x}_{t}^{i}|^{2}\right)\leq \frac{C}{n}
\]
where $C$ does not depend on $n$.
\end{pro}
\begin{rem}
In Proposition~\ref{mpmm2014pro1} and \ref{mpmm2014pro2}, we only state 1-dimensional case. In fact, the results
can be similarly generalized to multi-dimensional case.
\end{rem}
\noindent In \citet{jourdain2007nonlinear}, the coefficients are defined on $\mathds{R}\times\mathcal{P}(\mathds{R})$
$(\mathds{R}^{k}\times\mathcal{P}(\mathds{R}^{k}))$. We consider an extended case in which the coefficients
are defined on $[0,T]\times\Omega\times\mathds{R}\times\mathcal{P}(\mathds{R})$, depending on variation of
the time and sample. Proofs of Propositions~\ref{mpmm2014pro1} and~\ref{mpmm2014pro2} are similar to
\citet{jourdain2007nonlinear}. With the help of \citet{protter2004stochastic} on general stochastic
differential equations, we mimic proofs of~\citet{jourdain2007nonlinear} to prove our results.
For detailed proofs, please see the Appendix.

Now, we apply Proposition~\ref{mpmm2014pro1} and related discussion to analyze solvability of the
system~(\ref{mpmm2014eqn1}) and equation~(\ref{mpmm2014eqn2}). Since $b$ and $\sigma$ are Lipschitz
continuous with respect to $x$, it remains to verify they are also Lipschitz continuous with respect
to the Vasertein metric. Noticing that $b$, $\sigma$, $\psi$ and $\phi$ are all Lipschitz continuous,
we have
\[
\begin{array}{ll}
       & \left|b\left(t,\cdot,\int\psi(x)\mu(dx)\right)r(\alpha_{t}(\cdot))-
                 b\left(t,\cdot,\int\psi(y)\nu(dy)\right)r(\alpha_{t}(\cdot))\right|  \\
  \leq & C\left|\int\psi(x)d\mu(x)-\int\psi(y)d\nu(y)\right| \\
  \leq & Cd(\mu,\nu)
\end{array}
\]
and similarly for $\sigma r$. Hence, for given control,
under assumptions (A.1), (A.2) and (A.3), Proposition~\ref{mpmm2014pro1}
implies equation~(\ref{mpmm2014eqn2}) admits a unique strong solution. Replacing $\mu$ and $\nu$
by empirical measures in the above inequality together with equation~(\ref{mpmm2014eqn5}), we
have that the coefficients in the system~(\ref{mpmm2014eqn1}) satisfy the property of functional
Lipschitz. Thus, for given controls, the system~(\ref{mpmm2014eqn1}) admits a unique solution.

\section{Sufficient conditions for maximum principle}\label{mpmm2014scmp}

In this section, we develop sufficient conditions for maximum principle for Problem~\ref{mpmm2014prob1}.
Define the Hamiltonian
\[
\mathcal{H}(t,x,\mu,u,i,p,q)\triangleq h(t,x,\int\varphi d\mu,u)r(i)
        + b(t,x,\int\psi d\mu,u)r(i)p + \sigma(t,x,\int\phi d\mu, u)r(i)q
\]
For notational convenience, whenever $x$ is random variable associated with probability law $\mu$, we rewrite
the Hamiltonian as
\[
H(t,x,u,\alpha,p,q)\triangleq h(t,x,\e\varphi(x),u)r(\alpha)
        + b(t,x,\e\psi(x),u)r(\alpha)p + \sigma(t,x,\e\phi(x), u)r(\alpha)q.
\]
We also denote by $b_{x}$, $b_{y}$ and $b_{v}$ the derivative of $b$ with respect to the state
variable, the "expected value" and the control variable, respectively, and similarly for $\sigma$,
$h$, $g$ and $\psi$, $\phi$, $\varphi$, $\chi$. We shall use short-hand notation
$b(t)=b(t,x_{t},\e(\psi(x_{t})),u_{t})$ and similarly for other functions.
Let $\hat{u}_{t}$ be an equilibrium strategy to
Problem~\ref{mpmm2014prob1} and $\hat{x}_{t}$ be the associated state variable. We define
$\hat{\psi}(t)\triangleq \psi(\hat{x}_{t})$ and $\hat{b}(t)\triangleq b(t,\hat{x}_{t},\e\hat{\psi}(t),\hat{u}_{t})$
and similarly for the other functions and their derivatives.
Then, the adjoint equation is given by
\[
\left\{
\begin{array}{lll}
  dp(t) & = & -\left(\hat{b}_{x}(t)r(\alpha_{t})\hat{p}_{t}
                        + \hat{\sigma}_{x}(t)r(\alpha_{t})\hat{q}_{t} + \hat{h}_{x}(t)r(\alpha_{t})\right)dt \\
        &   & -\left(\e\left(\hat{b}_{y}(t)r(\alpha_{t})\hat{p}_{t}\right)\hat{\psi}_{x}(t) +
                     \e\left(\hat{\sigma}_{y}(t)r(\alpha_{t})\hat{q}_{t}\right)\hat{\phi}_{x}(t)
                          + \e\left(\hat{h}_{y}(t)r(\alpha_{t})\right)\hat{\varphi}_{x}(t)\right)dt  \\
        &   & + \hat{q}_{t}dw_{t} + s(t)d\tilde{\Phi}_{t} \\
  \hat{p}_{T} & = & \hat{g}_{x}(T)r(\alpha_{T}) + \e(\hat{g}_{y}(T)r(\alpha_{T}))\hat{\chi}_{x}(T),
\end{array}
\right.
\]
which is a backward stochastic differential equation (BSDE). For the existence and uniqueness of solutions
to BSDEs, see \citet{pardoux1990adapted}. For existence and uniqueness of solutions to BSDEs
driven by Markov chains, see \citet{cohen2010comparisons}.

We impose the following assumptions for sufficient conditions for maximum principle:
\begin{description}
\item[(A.4)] the function $g$ is convex in $(x,y)$.
\item[(A.5)] the Hamiltonian is convex in $(x,y,v)$.
\item[(A.6)] $\psi, \phi, \varphi, \chi$ are convex.
\item[(A.7)] the functions $b_{y}, \sigma_{y}, h_{y}, g_{y}$ are nonnegative.
\end{description}
\begin{thm}\label{mpmm2014thm1}
Assume the assumptions (A.1)--(A.7) are satisfied and let
$\hat{u}\in L_{\bar{\mathcal{F}}}^{2}(0,T;\mathds{U})$ with corresponding state process $\hat{x}_{t}$
and suppose there exists solutions $(\hat{p}_{t},\hat{q}_{t},\hat{s}_{t})$ to the adjoint equation satisfying
for all $u\in L_{\bar{\mathcal{F}}}^{2}(0,T;\mathds{U})$,
\eqb\label{mpmm2014eqn6}
\e\int_{0}^{T}|\sigma(t)r(\alpha_{t})\hat{p}_{t}|^{2}dt < \infty
\eqe
\eqb\label{mpmm2014eqn7}
\e\int_{0}^{T}|(\hat{x}_{t}-x_{t})\hat{q}_{t}|^{2}dt < \infty
\eqe
\eqb\label{mpmm2014eqn8}
\e\int_{0}^{T}|(\hat{x}_{t}-x_{t})\hat{s}_{t}^{*}Diag(m_{t})\hat{s}_{t}(\hat{x}_{t}-x_{t})|dt < \infty
\eqe
Then, if
\eqb\label{mpmm2014eqn9}
H(t,\hat{x}_{t},\hat{u}_{t},\alpha_{t},\hat{p}_{t},\hat{q}_{t})=
\inf_{v} H(t,\hat{x}_{t},v,\alpha_{t},\hat{p}_{t},\hat{q}_{t})
\eqe
for all $t\in [0,T]$, $\mathbb{P}$-a.s., $\hat{u}$ is an optimal strategy to Problem~\ref{mpmm2014pro1}.
\end{thm}
\Proof Let $H(t)\triangleq H(t,x_{t},u_{t},\alpha_{t},\hat{p}_{t},\hat{q}_{t})$ and
$\hat{H}(t)\triangleq H(t,\hat{x}_{t},\hat{u}_{t},\alpha_{t},\hat{p}_{t},\hat{q}_{t})$.
For any $u\in L_{\bar{\mathcal{F}}}^{2}(0,T;\mathds{U})$, we have
\[
J(\hat{u}) - J(u) = \e\int_{0}^{T}\left(\hat{h}(t)-h(t)\right)r(\alpha_{t})dt +
\e (\left(\hat{g}(T)-g(T)\right)r(\alpha_{T})).
\]

By the convexity of $g$ and $\chi$ as well as $g_{y}\geq 0$ and $r>0$, we obtain
\[
\begin{array}{lll}
  \e((\hat{g}-g)r) & \leq & \e\left(\hat{g}_{x}(T)r(\alpha_{T})(\hat{x}_{T}-x_{T})+
                          \hat{g}_{y}(T)r(\alpha_{T})\e(\hat{\chi}(T)-\chi(T))\right) \\
                 & \leq & \e\left(\hat{g}_{x}(T)r(\alpha_{T})(\hat{x}_{T}-x_{T})+
                          \hat{g}_{y}(T)r(\alpha_{T})\e(\hat{\chi}_{x}(T)(\hat{x}_{T}-x_{T}))\right) \\
                 & =    & \e(\hat{p}_{T}(\hat{x}_{T}-x_{T})).
\end{array}
\]
Apply Ir\^{o}'s formula to expand $\hat{p}(\hat{x}_{T}-x_{T})$ to get
\[
\begin{array}{lll}
  \hat{p}_{T}(\hat{x}_{T}-x_{T}) & = & \int_{0}^{T}(\hat{x}_{t}-x_{t})d\hat{p}_{t}
                                      + \int_{0}^{T}\hat{p}_{t}d(\hat{x}_{t}-x_{t})
                                      + [\hat{p}_{t}, \hat{x}_{t}-x_{t}](T)\\
                                 & = & \int_{0}^{T}(\hat{x}_{t}-x_{t})\left\{-\left(\hat{b}_{x}(t)r(\alpha_{t})\hat{p}_{t}
                       + \hat{\sigma}_{x}(t)r(\alpha_{t})\hat{q}_{t} + \hat{h}_{x}(t)r(\alpha_{t})\right)dt \right.\\
                                 &   & -\left(\e\left(\hat{b}_{y}(t)r(\alpha_{t})\hat{p}_{t}\right)\hat{\psi}_{x}(t) +
                     \e\left(\hat{\sigma}_{y}(t)r(\alpha_{t})\hat{q}_{t}\right)\hat{\phi}_{x}(t)
                          + \e\left(\hat{h}_{y}(t)r(\alpha_{t})\right)\hat{\varphi}_{x}(t)\right)dt  \\
                                 &   &  \left. + \hat{q}_{t}dw_{t} + s(t)d\tilde{\Phi}_{t}\right\}\\
                                 &   &   + \int_{0}^{T}\hat{p}_{t}(\hat{b}(t)-b(t))r(\alpha_{t})dt
                                        +  \int_{0}^{T}\hat{p}_{t}(\hat{\sigma}(t)-\sigma(t))r(\alpha_{t})dt\\
                                 &   &  + \int_{0}^{T}\hat{q}_{t}(\hat{\sigma}(t)-\sigma(t))r(\alpha_{t})dt .
\end{array}
\]
Due to the integrability condition~(\ref{mpmm2014eqn6}), (\ref{mpmm2014eqn7}) and~(\ref{mpmm2014eqn8}),
the Brownian motion and Markov chain martingale integrals in the above equation are square integrable
martingales which are null at the origin, we have
\[
\begin{array}{ll}
    & \e(\hat{p}_{T}(\hat{x}_{T}-x_{T})) \\
  = & -\e\int_{0}^{T}(\hat{x}_{t}-x_{t})(\hat{b}_{x}(t)r(\alpha_{t})\hat{p}_{t}
            +\e\left(\hat{b}_{y}(t)r(\alpha_{t})\hat{p}_{t}\right)\hat{\psi}_{x}(t) \\
    & + \hat{\sigma}_{x}(t)r(\alpha_{t})\hat{q}_{t}
      +\e\left(\hat{\sigma}_{y}(t)r(\alpha_{t})\hat{q}_{t}\right)\hat{\phi}_{x}(t)
      +\hat{h}_{x}(t)r(\alpha_{t})
      +\e\left(\hat{h}_{y}(t)r(\alpha_{t})\right)\hat{\varphi}_{x}(t))dt\\
    & +\e\int_{0}^{T}\hat{p}_{t}(\hat{b}(t)-b(t))r(\alpha_{t})dt
      + \int_{0}^{T}\hat{q}_{t}(\hat{\sigma}(t)-\sigma(t))r(\alpha_{t})dt
\end{array}
\]
Hence,
\[
\begin{array}{ll}
    & J(\hat{u})-J(u) \\
  \leq & \e\int_{0}^{T}\left(\hat{h}(t)-h(t)\right)r(\alpha_{t})dt + \e(\hat{p}_{T}(\hat{x}_{T}-x_{T})) \\
  =   & \e\int_{0}^{T}\left(\hat{H}(t)-H(t)\right)dt
         -\e\int_{0}^{T}\hat{p}_{t}\left(\hat{b}(t)-b(t)\right)r(\alpha_{t})dt
         -\e\int_{0}^{T}\hat{q}_{t}\left(\hat{\sigma}(t)-\sigma(t)\right)r(\alpha_{t})dt\\
      & + \e(\hat{p}_{T}(\hat{x}_{T}-x_{T})) \\
  =   &  \e\int_{0}^{T}\left(\hat{H}(t)-H(t)\right)dt
         -\e\int_{0}^{T}(\hat{x}_{t}-x_{t})(\hat{b}_{x}(t)r(\alpha_{t})\hat{p}_{t}
            +\e\left(\hat{b}_{y}(t)r(\alpha_{t})\hat{p}_{t}\right)\hat{\psi}_{x}(t) \\
      &  + \hat{\sigma}_{x}(t)r(\alpha_{t})\hat{q}_{t}
      +\e\left(\hat{\sigma}_{y}(t)r(\alpha_{t})\hat{q}_{t}\right)\hat{\phi}_{x}(t)
      +\hat{h}_{x}(t)r(\alpha_{t})
      +\e\left(\hat{h}_{y}(t)r(\alpha_{t})\right)\hat{\varphi}_{x}(t))dt.
\end{array}
\]

On the other hand, we differentiate the Hamiltonian and use the convexity of the functions to get
for all $t\in [0,T]$, $\mathbb{P}$-$a.s.$
\[
\begin{array}{l}
  \left(\hat{H}(t)-H(t)\right) \\
     \leq  \hat{H}_{x}(t)(\hat{x}_{t}-x_{t})
              + \hat{h}_{y}(t)r(\alpha_{t})\e(\hat{\varphi}(t)-\varphi(t))\\
           +\hat{b}_{y}(t)r(\alpha_{t})\e(\hat{\psi}(t)-\psi(t))\hat{p}_{t}
                 +\hat{\sigma}_{y}(t)r(\alpha_{t})\e(\hat{\phi}(t)-\phi(t))\hat{q}_{t}
                   + \hat{H}_{u}(t)(\hat{u}_{t}-u_{t})\\
     \leq  \hat{H}_{x}(t)(\hat{x}_{t}-x_{t})
               + \hat{h}_{y}(t)r(\alpha_{t})\e(\hat{\varphi}_{x}(t)(\hat{x}_{t}-x_{t})) \\
           +\hat{b}_{y}(t)r(\alpha_{t})\e(\hat{\psi}_{x}(t)(\hat{x}_{t}-x_{t}))\hat{p}_{t}
                 +\hat{\sigma}_{y}(t)r(\alpha_{t})\e(\hat{\phi}_{x}(t)(\hat{x}_{t}-x_{t}))\hat{q}_{t}
                   + \hat{H}_{u}(t)(\hat{u}_{t}-u_{t}) \\
     \leq   \hat{H}_{x}(t)(\hat{x}_{t}-x_{t})
                 + \hat{h}_{y}(t)r(\alpha_{t})\e(\hat{\varphi}_{x}(t)(\hat{x}_{t}-x_{t}))\\
           +\hat{b}_{y}(t)r(\alpha_{t})\e(\hat{\psi}_{x}(t)(\hat{x}_{t}-x_{t}))\hat{p}_{t}
               +\hat{\sigma}_{y}(t)r(\alpha_{t})\e(\hat{\phi}_{x}(t)(\hat{x}_{t}-x_{t}))\hat{q}_{t}
\end{array}
\]
where in the last inequality, we have used that $\hat{H}_{u}(t)(\hat{u}_{t}-u_{t})\leq 0$ due to
the minimization condition~(\ref{mpmm2014eqn9}). Therefore, $J(\hat{u})-J(u) \leq 0$. Thus,
$\hat{u}$ is an optimal strategy.
\qed




\section{Approximate Nash equilibrium}\label{mpmm2014ane}

We assume that the assumptions of Theorem~\ref{mpmm2014thm1} are satisfied,
therefore, Problem~\ref{mpmm2014prob1} has an optimal
strategy. In this section, we shall show that the optimal feedback strategy
$u_{t}^{i}\triangleq u(t,x_{t}^{i},\alpha_{t}^{i})$ solved from Problem~\ref{mpmm2014prob1}
is an $\varepsilon$-Nash equilibrium. In order to show that $(u^{1},\dots,u^{n})$ is an $\varepsilon$-Nash
equilibrium, we prove that for any $\varepsilon>0$, there exists $N>0$ such that whenever
$n>N$, the definition~\ref{mpmm2014defn1} is satisfied. We impose an additional assumption for
feedback control strategy: for $i=1,\dots,n$, $u^{i}$ satisfies Lipschitz condition. We write the
closed-loop equation and the associated Mckean-Vlasov equation as follows: for $i=1,\dots,n$,
\eqb\label{mpmm2014eqn10}
\begin{array}{ll}
  dx_{t}^{i,n} = & b\left(t, x_{t-}^{i,n}, \frac{1}{n}\sum\limits_{j=1}^{n}\psi(x_{t-}^{j,n}),
                          u(t,x_{t-}^{i,n},\alpha_{t-}^{i})\right)r\left(\alpha_{t-}^{i}\right)dt \\
                 & \hspace{20mm} +\sigma\left(t, x_{t-}^{i,n}, \frac{1}{n}\sum\limits_{j=1}^{n}\phi(x_{t-}^{j,n}),
                          u(t,x_{t-}^{i,n},\alpha_{t-}^{i})\right)r\left(\alpha_{t-}^{i}\right)dw_{t}^{i}.
\end{array}
\eqe
\eqb\label{mpmm2014eqn11}
    dx_{t}^{i}  =  b(t, x_{t-}^{i}, \e\psi(x_{t-}^{i}), u(t,x_{t-}^{i},\alpha_{t-}^{i}))r(\alpha_{t-}^{i})dt
                         +\sigma(t, x_{t-}^{i}, \e\phi(x_{t-}^{i}), u(t,x_{t-}^{i},\alpha_{t-}^{i}))r(\alpha_{t-}^{i})dw_{t}^{i}.
\eqe
Then, we have that $x^{i,n}$ can be approximated by $x^{i}$ as $n\to\infty$.
\begin{pro}\label{mpmm2014pro3}
As $n\to\infty$, we have that
$
\sup_{i\leq n}\e\left(\sup_{t\leq T}|x_{t}^{i,n}-x_{t}^{i}|^{2}\right)\to 0.
$
Moreover, if $b$ and $\sigma$ is linear in the third variable as in Proposition~\ref{mpmm2014pro2},
then,
$
\sup_{i\leq n}\e\left(\sup_{t\leq T}|x_{t}^{i,n}-x_{t}^{i}|^{2}\right)= O\left(\frac{1}{n}\right)
$
\end{pro}
For the proof, see Appendix 2. Now, let $\varepsilon_{n}^{1}\triangleq
\sup_{i\leq n}\e\left(\sup_{t\leq T}|x_{t}^{i,n}-x_{t}^{i}|^{2}\right)$. Proposition~\ref{mpmm2014pro3} implies
$\lim_{n\to\infty}\varepsilon_{n}\to 0$ and
$\mathcal{J}^{i}(u^{1}, \dots, u^{n}) = J^{i}(u^{i}) + O(\sqrt{\varepsilon_{n}})$, where $J^{i}(u^{i})$
formulated as Problem~\ref{mpmm2014prob1} is
the limiting optimization problem corresponding to $\mathcal{J}^{i}(u^{1}, \dots, u^{n})$. $\varepsilon_{n}$
will be determined.

\begin{thm}\label{mpmm2014thm2}
$(u^{1},\dots,u^{n})$ is an $\varepsilon$-Nash equilibrium of the cost $\mathcal{J}^{i}$ subject to
the system~(\ref{mpmm2014eqn1}), for $i=1,\dots,n$. That is,
for any fixed $1\leq i \leq n$, we have
\[
\mathcal{J}^{i}(u^{1}, \dots, u^{n}) \leq \mathcal{J}^{i}(u^{1}, \dots, u^{i-1}, v^{i}, u^{i+1}, \dots, u^{n})
+\sqrt{\varepsilon_{n}},
\]
when any alternative control
$v^{i}\in L_{\mathcal{F}^{i}}^{2}(0, T; \mathds{U})$ which is another Lipschitz feedback
is applied by the $i$th agent.
\end{thm}
\begin{rem}
If $b$, $\sigma$, $h$ and $g$ are linear in the third variable as in Proposition~\ref{mpmm2014pro2}, then
$\varepsilon_{n}$ will be specified as $\frac{1}{n}$.
\end{rem}
For the proof of Theorem~\ref{mpmm2014thm2}, see Appendix 2. We simply interpret the above theorem as follows.
If a given agent changes its control, it results in state process variations for other agents. These
variations and the initial control will affect the dynamics of that agent.


\section{Concluding remark}\label{mpmm2014conr}

We have proved sufficient stochastic maximum principle for Markov-modulated diffusion model.
We modulate the dynamics in a special way by multiplying coefficients by a positive function
of a Markov chain. It is possible that a analogy to more general coefficients involving a
Markov chain which satisfy suitable Lipchitz condition. On the other hand, a generalization
to Markov-modulated jump diffusion for maximum principle is also possible. On the other hand,
developing necessary conditions for maximum principle is also possible.

\newpage
\section*{Appendix 1: proofs for solutions of mean-field type SDEs}\label{mpmm2014psmf}

{\bf Proof of Proposition~\ref{mpmm2014pro1}}: Let $\mathds{D}$ denote the space of c\`{a}dl\`{a}g
functions from $[0,T]$ to $\mathds{R}$, $\mathcal{P}_{2}(\mathds{D})$ the space of probability measures
$\mathbb{Q}$ on $\mathds{D}$ such that
$\int_{\mathds{D}}\sup_{t\leq T}|\tilde{y}_{t}|^{2}\mathbb{Q}(d\tilde{y}) < \infty$.
For $\mathbb{P},\, \mathbb{Q}\in\mathcal{P}_{2}(\mathds{D})$, define the Vaserstein metric, for $t\in [0,T]$,
\[
\begin{array}{l}
     D_{t}(\mathbb{P},\mathbb{Q}) \\
    =\inf\left\{\left(\int_{\mathds{D}\times\mathds{D}}\sup_{s\leq t}
       |\tilde{y}_{s}-\tilde{w}_{s}|^{2}R(d\tilde{y},d\tilde{w})\right)^{1/2}:
       R\in\mathcal{P}(\mathds{D}\times\mathds{D})\,\,
       \text{with marginals}\,\, \mathbb{P}\,\, \text{and}\,\, \mathbb{Q}\right\}.
\end{array}
\]
Under the above metric, $\mathcal{P}_{2}(\mathds{D})$ is a complete space.

For any fixed $\mathbb{Q}\in\mathcal{P}_{2}(\mathds{D})$ with time-marginals $\{\mathbb{Q}_{t}: t\in [0,T]\}$,
we show
\eqb\label{mpmm2014a1}
\tilde{x}_{t}^{\mathbb{Q}}=\tilde{x}(0) +
\int_{0}^{t}\tilde{\sigma}(s,\tilde{x}_{s-}^{\mathbb{Q}},\mathbb{Q}_{s-})r(\tilde{\alpha}_{s-})dz_{s},\,\,
t\in [0,T]
\eqe
admits a unique solution, where
$\mathbb{Q}_{s-}=\mathbb{Q}\circ\tilde{y}_{s-}^{-1}$ is the weak limit of $\mathbb{Q}_{t}$
as $t\to s$ increasingly. Noticing that by Lebesgue's theorem, as $s\to t$ increasingly, the distance
\[
d(\mathbb{Q}_{t-},\mathbb{Q}_{s-})\leq\int_{\mathds{D}}|\tilde{y}_{t-}-\tilde{y}_{s-}|^{2}\mathbb{Q}(d\tilde{y})
\]
converges to 0. Similarly,
\[
d(\mathbb{Q}_{t},\mathbb{Q}_{s-})\leq\int_{\mathds{D}}|\tilde{y}_{t}-\tilde{y}_{s-}|^{2}\mathbb{Q}(d\tilde{y})
\]
converges to 0 as $s\to t$ decreasingly. Hence, we obtain that the mapping $t\in [0,T]\longrightarrow\mathbb{Q}_{t}$
is c\`{a}dl\`{a}d under the metric $d$ defined on $\mathcal{P}_{2}(\mathds{R})$. Thus, for fixed
$x\in \mathds{R}$, the mapping $t\in [0,T]\longrightarrow\tilde{\sigma}(t,x,\mathbb{Q}_{t})$ is c\`{a}dl\`{a}d.
On the other hand, $r(\cdot)$ is continuous function and $\tilde{\alpha}_{t-}$ is c\`{a}dl\`{a}d.
Hence, $\tilde{\sigma}(t,x,\mathbb{Q}_{t})r(\tilde{\alpha}_{t}(\omega))$ is c\`{a}dl\`{a}d. Then,
according to Theorem 6, p. 249, in \citet{protter2004stochastic},
equation~(\ref{mpmm2014a1}) admits a unique strong solution since
$\tilde{\sigma}(t,x,\mathbb{Q}_{t})r(\tilde{\alpha}_{t}(\omega))$ is random continuous.

Let $\Phi$ denote the mapping on $\mathcal{P}_{2}(\mathds{D})$ which associates the law of
$\tilde{x}^{\mathbb{Q}}$ with $\mathbb{Q}$. To use fixed point method, we shall verify that
$\Phi$ is a mapping from $\mathcal{P}_{2}(\mathds{D})$ to $\mathcal{P}_{2}(\mathds{D})$.
Indeed, for $K>0$, we set $\tau_{K}=\inf\{s\leq T: |\tilde{x}_{s}^{\mathbb{Q}}|\geq K\}$.
By Theorem 66, p.339 in \citet{protter2004stochastic} and Lipschitz property of
$\tilde{\sigma}(t,\cdot,\cdot)r(\tilde{\alpha}_{t}(\cdot))$, we have
\[
\begin{array}{l}
  \e\left(\sup_{s\leq t}\left|\tilde{x}_{s\wedge\tau_{K}}^{Q}\right|^{2}\right) \\
  \leq C\left(\e\left|\tilde{x}(0)\right|^{2}
          + \int_{0}^{t}\e\left(\boldsymbol{1}_{\{s\leq\tau_{K}\}}
            \left|\tilde{\sigma}(s,\tilde{x}_{s}^{\mathbb{Q}},\mathbb{Q}_{s})-\tilde{\sigma}(s,0,\delta_{0})\right|^{2}
            \left|r(\tilde{\alpha}_{s})\right|^{2}
              + \left|\tilde{\sigma}(s,0,\delta_{0})\right|^{2} \left|r(\tilde{\alpha}_{s})\right|^{2}\right)ds\right)\\
  \leq C\left(\e\left|\tilde{x}(0)\right|^{2}
          + \int_{0}^{t}\e\left(\boldsymbol{1}_{\{s\leq\tau_{K}\}}
            \left|\tilde{\sigma}(s,\tilde{x}_{s}^{\mathbb{Q}},\mathbb{Q}_{s})-\tilde{\sigma}(s,0,\delta_{0})\right|^{2}
              + \left|\tilde{\sigma}(s,0,\delta_{0})\right|^{2} \right)ds\right)\\
  \leq C\left(\e\left|\tilde{x}(0)\right|^{2}
          + \int_{0}^{t}\e\left(\sup_{u\leq s}\left|\tilde{x}_{u\wedge\tau_{K}}\right|^{2}\right)ds
            + t\int_{\mathds{D}}\sup_{t\leq T}\left|\tilde{y}_{t}\right|^{2}\mathbb{Q}(d\tilde{y})
              + \int_{0}^{t}\left|\tilde{\sigma}(s,0,\delta_{0})\right|^{2}ds \right).
\end{array}
\]
By Gronwall's Lemma, it follows that
\[
\e\left(\sup_{s\leq t}\left|\tilde{x}_{s\wedge\tau_{K}}^{Q}\right|^{2}\right)\leq
        C\left(\e\left|\tilde{x}(0)\right|^{2}
            + \int_{\mathds{D}}\sup_{t\leq T}\left|\tilde{y}_{t}\right|^{2}\mathbb{Q}(d\tilde{y})
              + \int_{0}^{t}\left|\tilde{\sigma}(s,0,\delta_{0})\right|^{2}ds \right),
\]
where $C$ does not depend on $K$. Let $K\to\infty$, by Fatou's Lemma,
\eqb\label{mpmm2014a3}
  \begin{array}{lll}
\int_{\mathds{D}}\sup_{t\leq T}\left|\tilde{y}_{t}\right|^{2}d\Phi(\mathbb{Q})(\tilde{y})
      & =    & \e\left(\sup_{s\leq t}\left|\tilde{x}_{s}^{Q}\right|^{2}\right) \\
      & \leq &  C\left(\e\left|\tilde{x}(0)\right|^{2}
            + \int_{\mathds{D}}\sup_{t\leq T}\left|\tilde{y}_{t}\right|^{2}\mathbb{Q}(d\tilde{y})
              + \int_{0}^{t}\left|\tilde{\sigma}(s,0,\delta_{0})\right|^{2}ds \right).
  \end{array}
\eqe
Hence, $\Phi$ is a mapping from $\mathcal{P}_{2}(\mathds{D})$ to $\mathcal{P}_{2}(\mathds{D})$.

Since  a process $\{\tilde{x}_{t}: t\in [0,T]\}$ such that
$\e\left(\sup_{t\leq T}\left|\tilde{x}_{t}\right|^{2}\right)<\infty$ solves
equation~(\ref{mpmm2014eqn3}) if and only if its law is a fixed point of $\Phi$. In the following,
we shall verify that $\Phi$ admits a unique fixed point. For
$\mathbb{P},\,\mathbb{Q}\in\mathcal{P}_{2}(\mathds{D})$, by a localization procedure similar to
the one used above, we have that
\[
\begin{array}{lll}
  \e\left(\sup_{s\leq t}\left|\tilde{x}_{s}^{\mathbb{P}}-\tilde{x}_{s}^{\mathbb{Q}}\right|^{2}\right)
   & \leq & C\int_{0}^{t}\e\left(\left|\tilde{\sigma}(s,\tilde{x}_{s}^{\mathbb{P}},\mathbb{P}_{s})
      -\tilde{\sigma}(s,\tilde{x}_{s}^{\mathbb{Q}},\mathbb{Q}_{s})\right|^{2}\left|r(\tilde{\alpha}_{s})\right|^{2}\right)ds \\
   & \leq & C\int_{0}^{t}\e\left(\left|\tilde{\sigma}(s,\tilde{x}_{s}^{\mathbb{P}},\mathbb{P}_{s})
                            -\tilde{\sigma}(s,\tilde{x}_{s}^{\mathbb{\mathbb{Q}}},\mathbb{Q}_{s})\right|^{2}\right)ds\\
   & \leq & C\int_{0}^{t}\left(\e\left(\sup_{u\leq s}
                    \left|\tilde{x}_{u}^{\mathbb{P}}-\tilde{x}_{u}^{\mathbb{Q}}\right|^{2}\right)
            + d^{2}(\mathbb{P}_{s},\mathbb{Q}_{s})\right)ds
\end{array}
\]
From Gronwall's Lemma, we obtain that $\forall\in [0,T]$,
\[
\e\left(\sup_{s\leq t}\left|\tilde{x}_{s}^{\mathbb{P}}-\tilde{x}_{s}^{\mathbb{Q}}\right|^{2}\right)
\leq C\int_{0}^{t}d^{2}(\mathbb{P}_{s},\mathbb{Q}_{s})ds.
\]
It is noted that $D_{t}^{2}\left(\Phi(\mathbb{P}),\Phi(\mathbb{Q})\right)\leq
\e\left(\sup_{s\leq t}\left|\tilde{x}_{s}^{\mathbb{P}}-\tilde{x}_{s}^{\mathbb{Q}}\right|^{2}\right)$
and $d(\mathbb{P}_{s},\mathbb{Q}_{s})\leq D_{s}(\mathbb{P},\mathbb{Q})$, we have that
$\forall t\in [0,T]$,
\[
D_{t}^{2}\left(\Phi(\mathbb{P}),\Phi(\mathbb{Q})\right)\leq
\int_{0}^{t}D_{s}^{2}\left(\mathbb{P},\mathbb{Q}\right)ds.
\]
Iterating this inequality and denoting by $\Phi^{n}$ the $n$-fold composition of $\Phi$, we obtain that
$n=1,2,\dots$,
\[
D_{T}^{2}\left(\Phi^{n}(\mathbb{P}),\Phi^{n}(\mathbb{Q})\right)\leq
C^{n}\int_{0}^{T}\frac{(T-s)^{n-1}}{(n-1)!}D_{s}^{2}\left(\mathbb{P},\mathbb{Q}\right)ds
\leq \frac{C^{n}T^{n}}{n!}D_{T}^{2}\left(\mathbb{P},\mathbb{Q}\right).
\]
Hence, for sufficiently large $n$, $\Phi^{n}$ is a contraction, therefore, $\Phi$ admits a
unique fixed point.
\qed

\noindent{\bf Proof of Proposition~\ref{mpmm2014pro2}:} Let
$\mathbb{P}^{n}=\frac{1}{n}\sum_{j=1}^{n}\delta_{\tilde{x}^{j}}$ be the empirical measure of the
independent nonlinear process
\[
\left\{
  \begin{array}{lll}
    d\tilde{x}_{t}^{i} & = & \tilde{\sigma}(t, \tilde{x}_{t-}^{i}, \p_{t}))
                              r(\tilde{\alpha}_{t-}^{i})dz_{t}, \quad t\in [0,T]\\
    \tilde{x}_{0}^{i}  & = & \tilde{x}_{0}^{i},\quad \forall t\in [0,T],\, \p_{t}\,\, \text{denotes the probability
                        distribution of}\,\, \tilde{x}_{t}^{i}.
  \end{array}
\right.
\]
By a localization procedure similar to the one used in the proof of Proposition~\ref{mpmm2014pro1},
we have
\eqb\label{mpmm2014eqna2}
\begin{array}{lll}
  \e\left(\sup_{s\leq t}\left|\tilde{x}_{s}^{i,n}-\tilde{x}_{s}^{i}\right|^{2}\right)
   & \leq & C\int_{0}^{t}\e\left(\left|\tilde{\sigma}(s,\tilde{x}_{s}^{i,n},\mu_{s}^{n})
   -\tilde{\sigma}(s,\tilde{x}_{s}^{i},\mathbb{P}_{s}^{n})\right|^{2}\left|r(\tilde{\alpha}_{s}^{i})\right|^{2}\right)ds \\
   &      & + \,\, C\int_{0}^{t}\e\left(\left|\tilde{\sigma}(s,\tilde{x}_{s}^{i},\mathbb{P}_{s}^{n})
                            -\tilde{\sigma}(s,\tilde{x}_{s}^{i},\mathbb{P}_{s})\right|^{2}\left|r(\tilde{\alpha}_{s}^{i})\right|^{2}\right)ds\\
   & \leq & C\int_{0}^{t}\e\left(\left|\tilde{\sigma}(s,\tilde{x}_{s}^{i,n},\mu_{s}^{n})
   -\tilde{\sigma}(s,\tilde{x}_{s}^{i},\mathbb{P}_{s}^{n})\right|^{2}\right)ds \\
   &      & + \,\, C\int_{0}^{t}\e\left(\left|\tilde{\sigma}(s,\tilde{x}_{s}^{i},\mathbb{P}_{s}^{n})
                            -\tilde{\sigma}(s,\tilde{x}_{s}^{i},\mathbb{P}_{s})\right|^{2}\right)ds
\end{array}
\eqe
Due to the Lipschitz property of $\tilde{\sigma}$, equation~(\ref{mpmm2014eqn5}) and
exchangeability of the couples $(\tilde{x}^{i},\tilde{x}^{i,n}), i=1,\dots,n$, the first term
of the right in the above inequality is less than
$C\int_{0}^{t}\e\left(\sup_{u\leq s}\left|\tilde{x}_{u}^{i,n}-\tilde{x}_{u}^{i}\right|^{2}\right)ds$.
By Gronwall's Lemma and Lipschitz assumption on $\tilde{\sigma}$, we have
\[
\e\left(\sup_{s\leq t}\left|\tilde{x}_{s}^{i,n}-\tilde{x}_{s}^{i}\right|^{2}\right)
\leq  C\int_{0}^{t}\e\left(\left|\tilde{\sigma}(s,\tilde{x}_{s}^{i},\mathbb{P}_{s}^{n})
                            -\tilde{\sigma}(s,\tilde{x}_{s}^{i},\mathbb{P}_{s})\right|^{2}\right)ds
\leq  C\int_{0}^{t}\e(d^{2}(\mathbb{P}_{s}^{n},\mathbb{P}_{s}))ds,
\]
From Lemma 4 in \citet{jourdain2007nonlinear}, the upper bounds of the
second order moments in Proposition~\ref{mpmm2014pro1},
we yield the first assertion.

Moreover, if $\tilde{\sigma}(t,\tilde{x},\nu)=\int_{\mathds{R}}\eta(t,\tilde{x},\tilde{y})\nu(d\tilde{y})$, we have
$\e\left(\left|\tilde{\sigma}(s,\tilde{x}_{s}^{i},\mathbb{P}_{s}^{n})
                            -\tilde{\sigma}(s,\tilde{x}_{s}^{i},\mathbb{P}_{s})\right|^{2}\right)$
is equal to
\[
\frac{1}{n}\sum_{j,l=1}^{n}\e\left(\left[\eta(s,\tilde{x}_{s}^{i},\tilde{x}_{s}^{j})-
\int_{\mathds{R}}\eta(s,\tilde{x}_{s}^{i},\tilde{y})\mathbb{P}_{s}(d\tilde{y})\right]
\left[\eta(s,\tilde{x}_{s}^{i},\tilde{x}_{s}^{l})-
\int_{\mathds{R}}\eta(s,\tilde{x}_{s}^{i},\tilde{y})\mathbb{P}_{s}(d\tilde{y})\right]\right).
\]
By the independence of the random variables $\tilde{x}_{s}^{1},\dots,\tilde{x}_{s}^{n}$ with common law
$\mathbb{P}_{s}$, the expectation in the above summation vanishes as long as $j\neq l$.
As a consequence, the result follows. \qed

\section*{Appendix 2: proofs for Proposition~\ref{mpmm2014pro3}}\label{mpmm2014pp}

{\bf Proof of Proposition~\ref{mpmm2014pro3}:} Let
$\tilde{\sigma}(t,x_{t}^{i,n},\mu_{t}^{n},u_{t}^{i,n})\triangleq \left(b\left(t, x_{t}^{i,n},
\frac{1}{n}\sum\limits_{j=1}^{n}\psi(x_{t}^{j,n}),u(t,x_{t}^{i,n},\alpha_{t}^{i})\right)\right.,\\
\left.\sigma\left(t, x_{t}^{i,n}, \frac{1}{n}\sum\limits_{j=1}^{n}\phi(x_{t}^{j,n}),
                          u(t,x_{t}^{i,n},\alpha_{t}^{i})\right)\right)$,
$z_{t}^{i}=(t, w_{t}^{i})^{*}$,
$\tilde{\sigma}(t,x_{t}^{i},\mu_{t},u_{t}^{i})\triangleq \left(b\left(t, x_{t}^{i},
\e\psi(x_{t}^{i}),u(t,x_{t}^{i},\alpha_{t}^{i})\right)\right.,\\
\left.\sigma\left(t, x_{t}^{i}, \e\phi(x_{t}^{i}),
                          u(t,x_{t}^{i},\alpha_{t}^{i})\right)\right)$, where
$\mu_{t}^{n}=\frac{1}{n}\sum\limits_{j=1}^{n}\delta_{x_{t}^{j,n}}$,
$\mu_{t}$ is the marginal distribution of $x_{t}^{i}$, $u_{t}^{i,n}=u(t,x_{t}^{i,n},\alpha_{t}^{i})$
and $u_{t}^{i}=u(t,x_{t}^{i},\alpha_{t}^{i})$. Let
$\nu_{t}^{n}\triangleq \frac{1}{n}\sum\limits_{j=1}^{n}\delta_{x_{t}^{j}}$. Then, we can rewrite equations
(\ref{mpmm2014eqn10}) and (\ref{mpmm2014eqn11}) as follows:
\[
  dx_{t}^{i,n} = \tilde{\sigma}(t,x_{t}^{i,n},\mu_{t}^{n},u_{t}^{i,n})r(\alpha_{t}^{i})dz_{t}^{i}.
\]
\[
    dx_{t}^{i}  =  \tilde{\sigma}(t,x_{t}^{i},\mu_{t},u_{t}^{i})r(\alpha_{t}^{i})dz_{t}^{i}.
\]
Hence, we have
\[
\begin{array}{lll}
  x_{t}^{i,n}-x_{t}^{i} & = & \int_{0}^{t} \left[\tilde{\sigma}(s,x_{s}^{i,n},\mu_{s}^{n},u_{s}^{i,n})-
                    \tilde{\sigma}(t,x_{t}^{i},\mu_{t},u_{t}^{i})\right]r(\alpha_{t}^{i}) dz_{t}^{i} \\
    & = &  \int_{0}^{t} \left[\tilde{\sigma}(s,x_{s}^{i,n},\mu_{s}^{n},u_{s}^{i,n})-
                    \tilde{\sigma}(t,x_{t}^{i},\nu_{t}^{n},u_{t}^{i})+
                         \tilde{\sigma}(t,x_{t}^{i},\nu_{t}^{n},u_{t}^{i})-
                            \tilde{\sigma}(t,x_{t}^{i},\mu_{t},u_{t}^{i})\right]r(\alpha_{t}^{i}) dz_{t}^{i}.
\end{array}
\]
It is noted that $b$ and $\sigma$ are differentiable with respect to $(x,y,v)$, and thus satisfy Lipschitz condition.
Since $u^{i,n}$ and $u^{i}$ satisfy Lipschitz condition, by a similar argument to
the proof of Proposition~\ref{mpmm2014pro2}, we yield Proposition~\ref{mpmm2014pro3}.\qed

\noindent{\bf Proof of Theorem~\ref{mpmm2014thm2}:} Due to symmetry of index $i$, we only need
to consider a control strategy for the first agent. We first analyze the running cost, and then
the terminal cost by a similar procedure. We write the system with changed control variable for
the first agent as follows:
\[
  d\hat{x}_{t}^{1,n} = \tilde{\sigma}(t,\hat{x}_{t}^{1,n},\hat{\mu}_{t}^{n},\hat{v}_{t}^{1,n})r(\alpha_{t}^{1})dz_{t}^{1},
\]
\[
  d\hat{x}_{t}^{i,n} = \tilde{\sigma}(t,\hat{x}_{t}^{i,n},\hat{\mu}_{t}^{n},\hat{u}_{t}^{i,n})r(\alpha_{t}^{i})dz_{t}^{i}, \quad i=2,\dots,n,
\]
where $\hat{\mu}_{t}^{n}$ and $\hat{v}_{t}^{1,n},\hat{u}_{t}^{2,n},\dots,\hat{u}_{t}^{n,n}$
are defined along the same line with $\mu_{t}^{n}$ and $u_{t}^{i,n}$, for $i=1,\dots,n$,
in the proof of Proposition~\ref{mpmm2014pro3}.

For $i\neq 1$, we have the following estimate
\[
\sup_{2\leq j \leq n}\e\left(\sup_{s\leq T}|x_{s}^{i,n}-\hat{x}_{s}^{i,n}|^{2}\right)\leq \varepsilon_{n}^{2},\quad
\text{for sufficiently large}\,\, n.
\]
Indeed, the above estimate can be verified by Gronwall's lemma together with euqtion~(\ref{mpmm2014eqn5}).
The expectation in the above equation is less than
$C\int_{0}^{T}\e\left(\sup_{u\leq s}\left(\left|x_{u}^{i,n}-\hat{x}_{u}^{i,n}\right|^{2}
+\frac{1}{n}\left|x_{u}^{1,n}-\hat{x}_{u}^{1,n}\right|^{2}\right)\right)ds$
by treating $x_{u}^{1,n}$ and $\hat{x}_{u}^{1,n}$ as additional quantities and applying
equation~(\ref{mpmm2014eqn5}). Finiteness of
$\left|\tilde{x}_{u}^{1,n}-\hat{x}_{u}^{1,n}\right|^{2}$ can be analyzed as same as equation~(\ref{mpmm2014a3}).
Then, the above inequality follows from Gronewall's lemma. By the same procedure as in
the proof of Proposition~\ref{mpmm2014pro3} and treating $x_{u}^{1}$ and $\hat{x}_{u}^{1,n}$
as additional quantities, we also have the following estimate
\[
\sup_{2\leq j \leq n}\e\left(\sup_{s\leq T}|x_{s}^{i}-\hat{x}_{s}^{i,n}|^{2}\right)\leq \varepsilon_{n}^{3},\quad
\text{for sufficiently large}\,\, n.
\]
We construct a new equation
\[
  d\bar{x}_{t}^{1,n} = \tilde{\sigma}(t,\bar{x}_{t}^{1,n},\nu_{t}^{n},\bar{v}_{t}^{1})r(\alpha_{t}^{1})dz_{t}^{1}.
\]
Then, we have the following estimate
\[
\e\left(\sup_{s\leq T}|\bar{x}_{s}^{1,n}-\hat{x}_{s}^{1,n}|^{2}\right)\leq \varepsilon_{n}^{4},\quad
\text{for sufficiently large}\,\, n.
\]
Now, we define the equation corresponding
\[
  d\hat{x}_{t}^{1} = \tilde{\sigma}(t,\hat{x}_{t}^{1},\mu_{t},\hat{v}_{t}^{1})r(\alpha_{t}^{1})dz_{t}^{1}.
\]
Then, by a similar argument to the proof of Proposition~\ref{mpmm2014pro3}, we obtain
\[
\e\left(\sup_{s\leq T}|\hat{x}_{s}^{1}-\bar{x}_{s}^{1,n}|^{2}\right)\leq \varepsilon_{n}^{5},\quad
\text{for sufficiently large}\,\, n.
\]

Let $\varepsilon_{n}=\max\{\varepsilon_{n}^{1},\dots,\varepsilon_{n}^{5}\}$ and
$\bar{h}(t,x_{t}^{i,n},\mu_{t}^{n},u_{t}^{i,n})\triangleq
h\left(t, x_{t}^{i,n},
\frac{1}{n}\sum\limits_{j=1}^{n}\psi(x_{t}^{j,n}),u(t,x_{t}^{i,n},\alpha_{t}^{i})\right)$.
Based on above estimates, we obtain
\[
\begin{array}{ll}
       & \e\int_{0}^{T}\bar{h}(t,\hat{x}_{t}^{1,n},\hat{\mu}_{t}^{n},\hat{v}_{t}^{1,n})dt\\
  \geq & \e\int_{0}^{T}\bar{h}(t,\hat{x}_{t}^{1,n},\nu_{t}^{n},\hat{v}_{t}^{1,n})dt - O(\sqrt{\varepsilon_{n}}) \\
  \geq & \e\int_{0}^{T}\bar{h}(t,\hat{x}_{t}^{1},\nu_{t}^{n},\hat{v}_{t}^{1})dt - O(\sqrt{\varepsilon_{n}}) \\
  \geq & \e\int_{0}^{T}\bar{h}(t,\hat{x}_{t}^{1},\mu_{t},\hat{v}_{t}^{1})dt - O(\sqrt{\varepsilon_{n}}) \\
  \geq & \e\int_{0}^{T}\bar{h}(t,x_{t}^{1},\mu_{t},u_{t}^{1})dt - O(\sqrt{\varepsilon_{n}}),
\end{array}
\]
where last inequality results from the optimality assumption. Similarly, we can analyze the terminal cost. Hence, we get
\[
\begin{array}{lll}
  \mathcal{J}^{}(v^{1}, u^{2}, \dots, u^{n}) & = & J^{1}(v^{1}) - O(\sqrt{\varepsilon_{n}})\\
                                                                       & \geq & J^{1}(u^{1}) - O(\sqrt{\varepsilon_{n}})\\
                                                                       & = & \mathcal{J}^{1}(u^{1}, \dots, u^{n})- O(\sqrt{\varepsilon_{n}}),
\end{array}
\]
where last equality follows from Proposition~\ref{mpmm2014pro3}.
\qed
\bibliographystyle{agsm}       
\bibliography{BibTex-SOC-Tai}   

\end{document}